\documentclass[oneside,11pt]{article}
\usepackage{amscd, amssymb, amsmath, amsthm, amsfonts}
\usepackage{latexsym, graphics, graphicx,psfrag}
\usepackage[all]{xy}

\usepackage{graphicx}
\begin{document}
\title{Additivity of Heegaard genera of bounded surface sums}

\author {Ruifeng Qiu\footnote{Supported by a grant of NSFC (No.
10625102)},  Shicheng Wang\footnote{Supported by a grant of NSFC
(No. 10631060)}, Mingxing Zhang}
\date{}
\maketitle
\begin{abstract}

Let $M$ be a surface sum of 3-manifolds $M_1$ and $M_2$ along a
bounded connected surface $F$ and $\partial_{i}$ be the component of
$\partial M_{i}$ containing $F$. If $M_i$ has a high distance
Heegaard splitting, then any minimal Heegaard splitting of $M$ is
the amalgamation of those of $M^{1}, M^{2}$ and $M^{*}$, where
$M^{i}=M_{i}\setminus
\partial_{i}\times I$, and $M^{*}=\partial_{1}\times I\cup_{F}
\partial_{2}\times I$. Furthermore, once both
$\partial_i\setminus F$ are connected, then
$g(M)=Min\bigl\{g(M_{1})+g(M_{2}), \alpha\bigr\}$, where $\alpha=
g(M_{1})+g(M_{2})+1/2(2\chi(F)+2-\chi(\partial _{1})-\chi(\partial
_{2}))-Max\bigl\{g(\partial _{1}), g(\partial _{2})\bigl\}$; in
particular $g(M)=g(M_{1})+g(M_{2})$ if and only if $\chi(F)\geq
1/2Max\bigl\{\chi(\partial_{1}), \chi(\partial_{2})\bigr\}.$

The proofs rely on Scharlemann-Tomova's theorem.
\end{abstract}



\section{Introduction}

 \ \ \ \ \ All surfaces and 3-manifolds in this paper are assumed to be compact and
orientable.

Let $F$ be either a properly embedded  surface in a 3-manifold $M$
or a sub-surface of $\partial M$. If there is an essential simple
closed curve on $F$ which bounds a disk in $M$ or $F$ is a 2-sphere
which bounds a 3-ball in $M$, then we say $F$ is compressible;
otherwise, $F$ is said to be incompressible. If $F$ is an
incompressible surface not parallel to $\partial M$, then $F$ is
said to be essential. If $M$ contains an essential 2-sphere, then
$M$ is said to be reducible. If $\partial M$ is compressible, then
$M$ is said to be $\partial$-reducible.

Let $M$ be a  3-manifold. If there is a closed surface $S$ which
cuts $M$ into two compression bodies $V$ and $W$ with
$S=\partial_{+} W=\partial_{+} V$, then we say $M$ has a Heegaard
splitting, denoted by $M=V\cup_{S} W$; and $S$ is called a Heegaard
surface of $M$. Moreover, if the genus $g(S)$ of $S$ is minimal
among all the Heegaard surfaces of $M$, then $g(S)$ is called the
genus of $M$, denoted by $g(M)$.

Now let $M$ be a 3-manifold, and $F$ be a compact surface in $M$
which cuts $M$ into two 3-manifolds $M_{1}$ and $M_{2}$. Then $M$ is
called a surface sum of $M_{1}$ and $M_{2}$ along $F$, denoted by
$M=M_{1}\cup_{F} M_{2}$. Note that $F\subset\partial M_{i}$ for
$i=1, 2$.  A central topic in Heegaard splitting is to address
relations between $g(M_1)$, $g(M_2)$ and $g(M)$.

Suppose first that $F$ is a closed surface.  Let
$M_{i}=V_{i}\cup_{S_{i}} W_{i}$ be a Heegaard splitting for $i= 1,
2$. Then $M$ has a natural Heegaard splitting called the
amalgamation of $V_{1}\cup_{S_{1}} W_{1}$ and $V_{2}\cup_{S_{2}}
W_{2}$. From this view, $g(M)\leq g(M_{1})+g(M_{2})-g(F)$. If $F$ is
a 2-sphere, the so-called Haken's lemma claimed
$g(M)=g(M_{1})+g(M_{2})$. For $g(F)>0$, there are some examples to
show that it is possible that $g(M)\leq g(M_{1})+g(M_{2})-g(F)-n$
for any given $n>0$,  see [7] and [18]. Philosophically, in such
examples neither the gluing between $M_1$ and $M_2$ along $F$ nor
the Heegaard splitting of $M_i$ are complicated.

Under various different conditions describing  the complicated
gluing maps, the equality $g(M)=g(M_{1})+g(M_{2})-g(F)$ was proved,
see [1], [9], [10] and [19].  By invoke results of Hartshorn [3],
Scharlemann [13] and Scharlemann and Tomova [16], it is just proved
in [6] that $g(M)=g(M_{1})+g(M_{2})-g(F)$ if  $M_1$ and $M_2$ have
high  distance Heegaard splittings, where the distance of a Heegaard
splitting was introduced by Hempel [4].

Suppose that $F$ is a bounded surface. Then  it is easy to see
$g(M)\leq g(M_{1})+g(M_{2})$ (see Lemma 2.1). By the disk version of
Haken's lemma, $g(M)=g(M_{1})+g(M_{2})$ if $F$ is a disk. If $F$ is
an annulus, various results about if $g(M)=g(M_{1})+g(M_{2})$ hold
or not have been given, see [5], [8], [11] and [12].

In this paper we will address the additivity of Heegaard genus of
surface sum of 3-manifolds along a bounded surface $F$ with
$\chi(F)< 0$, which  seems not touched before.

We first fix some notions.  Suppose $P$ (resp. $H$) is a properly
embedded surface (resp. an embedded 3-manifold) in a 3-manifold $M$.
We use $M\setminus P$ (resp. $M\setminus H$) to denote the resulting
manifold obtained by splitting $M$ along $P$ (resp. removing
$\text{int} H$, the interior of $H$).

Let $M=M_{1}\cup_{F} M_{2}$, $\partial_{i}$ be the component of
$\partial M_{i}$ containing $F$, and $\partial_{i}\times [0,1]$ be a
regular neighborhood of $\partial_{i}$ in $M_{i}$ with $\partial
_{i}=\partial _{i}\times\bigl\{0\bigr\}$. We denote by $P^{i}$ the
surface $\partial_{i}\times\bigl\{1\bigr\}$. Let
$M^{i}=M_{i}\setminus \partial_{i}\times [0,1]$ for $i=1, 2$,  and
$M^{*}=\partial_{1}\times [0,1]\cup_{F}
\partial_{2}\times [0,1]$. Then
$M=M^{1}\cup_{P^{1}}M^{*}\cup_{P^{2}}M^{2}$.  \vskip 3mm

{\bf Theorem 1.} \ Let $M$ be  a surface sum of 3-manifolds $M_1$
and $M_2$ along a bounded surface $F$, and $\partial_{i}$ be the
component of $\partial M_{i}$ containing $F$. If $M_{i}$ has a
Heegaard splitting $V_{i}\cup_{S_{i}} W_{i}$ with $d(S_{i})>
2(g(S_{1})+g(S_{2}))$, $i=1,2$,  then any minimal Heegaard splitting
of $M$ is the amalgamation of Heegaard splittings of $M^{1}$,
$M^{2}$ and $M^{*}$ along $\partial_{1}$ and $\partial_{2}$. \vskip
3mm

The proof of Theorem 1  invokes full energy of Scharlemann-Tomova's
deep result (Lemma 2.5).

{\bf Theorem 2.} \ Under the assumptions of Theorem 1, if
$\partial_i \setminus F$ is connected for $i=1,2$,   then
$g(M)=Min\bigl\{g(M_{1})+g(M_{2}), \alpha\bigr\}$, where  $$\alpha=
g(M_{1})+g(M_{2})+1/2(2\chi(F)+2-\chi(\partial _{1})-\chi(\partial
_{2}))-Max\bigl\{g(\partial _{1}), g(\partial _{2})\bigl\}\qquad
(1.1).$$  Furthermore $g(M)=g(M_{1})+g(M_{2})$ if and only if
$\chi(F)\geq 1/2Max\bigl\{\chi(\partial_{1}),
\chi(\partial_{2})\bigr\}.$ \vskip 3mm

{\bf Corollary 3.} \ Under the assumptions of Theorem 1, if $F$ is
annulus, then $g(M)=g(M_{1})+g(M_{2})$.\vskip 3mm

{\bf Remark \ 4.} It is remarkable to compare Theorem 2 with the
main result in [6] for closed surface case.  In Theorem 2

(1) once $\chi(F)\geq 1/2Max\bigl\{\chi(\partial_{1}),
\chi(\partial_{2})\bigr\}$, $\chi(F)$, therefore $g(F)$, itself
plays no role in the result.

(2) once $\chi(F)< 1/2Min\bigl\{\chi(\partial_{1}),
\chi(\partial_{2})\bigr\}$, the contribution of $\chi(F)$ to $g(M)$
is non-trivial linear function. A particular case when
$g(\partial_{1})=g(\partial_{2})$ and $\chi(F)=\chi(\partial_{1})$,
then $g(M)=g(M_{1})+g(M_{2})-g(\partial_{1})+1$.

{\bf Acknowledgements.} The authors would like to thank Tsuyoshi
Kobayashi and Qing Zhou for helpful discussions on this paper.
\section{Distance of Heegaard splitting}

Weakly incompressible surface in 3-manifolds was introduced in [16]:
Let $P$ be a separating connected closed surface in 3-manifold $M$
which cuts $M$ into two 3-manifolds $M_{1}$ and $M_{2}$. Then $P$ is
said to be bicompressible if $P$ is compressible in both $M_{1}$ and
$M_{2}$. $P$ is strongly compressible if there are compressing disks
for $P$ in $M_{1}$ and $M_{2}$ which have disjoint boundaries in
$P$; otherwise $P$ is weakly incompressible.

Let $M$ be a compact orientable 3-manifold, and $M=V\cup_{S} W$ be a
Heegaard splitting. If there are essential disks $B\subset V$ and
$D\subset W$ such that $\partial B=\partial D$ ($\partial
B\cap\partial D=\emptyset$), then $V\cup_{S} W$ is said to be
reducible (weakly reducible). Otherwise, it is said to be
irreducible (strongly irreducible), see [2]. It is easy to see that
a strongly irreducible Heegaard surface is weakly incompressible.

Now let $P$ be a bicompressible surface in an irreducible 3-manifold
$M$. By maximally compressing $P$ in both sides of $P$ and deleting
the possible 2-sphere components, we get a surface sum structure of
$M$ as follow:
$$M=N_1\cup_{F^{P}_{1}}H^{P}_{1}\cup_{P}H^{P}_{2}\cup_{F^{P}_{2}}N_2,$$
where $H^{P}_{i}$ is a compression body with $\partial_{+}
H^{P}_{i}=P$, and  $\partial_{-} H^{P}_{i}=F^{P}_{i}$ is a
collection (may be empty) of incompressible closed surfaces of genus
at least one in $N_{i}$, $i=1,2$. Note that, if $F^{P}_{i}$ is
empty, then $H^{P}_{i}$ is a handlebody and $N_{i}$ is empty. It is
easy to see that if $M$ has boundary, then $F^{P}_{1}$ and
$F^{P}_{2}$ can not be both empty. Moreover if $P$ is weakly
incompressible, then the Heegaard splitting
$H^{P}_{1}\cup_{P}H^{P}_{2}$ is strongly irreducible.

Two weakly incompressible surfaces $P$ and $Q$ are said to be
well-separated in $M$ if $H^{P}_{1}\cup_{P}H^{P}_{2}$ is disjoint
from $H^{Q}_{1}\cup_{P}H^{Q}_{2}$ by isotopy. \vskip 2mm

{\bf Lemma 2.1.} \ Suppose $P$ is a weakly incompressible surface in
$M$. Then each component of $F^{P}_{i}$ is incompressible in $M$ for
$i=1, 2$.

{\bf Proof.} \ By the definition, each component of $F^{P}_{i}$ is
incompressible in $N_{i}$ for $i=1,2$. Since $P$ is a bicompressible
but weakly incompressible, $H^{P}_{1}\cup H^{P}_{2}$ is a
non-trivial strongly Heegaard splitting. By the disc version of
Haken's Lemma, each component of $F^{P}_{i}$ is incompressible in
$H^{P}_{1}\cup H^{P}_{2}$. Q.E.D. \vskip 2mm

{\bf Lemma 2.2.}  \ Let $S$ be a strongly irreducible Heegaard
surface of a 3-manifold $M$, then $M\setminus (H_{1}^{S}\cup_{S}
H_{2}^{S})$ is homeomorphic to $\partial M\times I$. Fruthermore, if
$P$ is a weakly incompressible surface in $M$, then either
$H^{P}_{1}\cup_{P} H^{P}_{2}\subset\partial \times I$ is
homeomorphic to $\partial\times I$ for one component $\partial$ of
$\partial M$, or, $S$ and $P$ are not well-separated.

{\bf Proof.} \ Now $M=V\cup_{S} W$, where $V$ and $W$ are
compression bodies. We may assume that $H_{1}^{S}\subset V$ and
$H_{2}^{P}\subset W$. Now there are essential disks $D_{1}, \ldots,
D_{n}$ in $V$ such that each component of $V \setminus
\cup_{i=1}^{n} D\times [0,1]$ is an I-bundle of a closed surface.
Since $\partial_{-} H_{1}^{S}$ is incompressible in $M$, hence in
$V$, each $D_{i}$ can be isotoped to be disjoint from $\partial_{-}
H_{1}^{S}$. Hence $\partial_{-} H_{1}^{S}\subset V \setminus
\cup_{i=1}^{n} D\times [0,1]$. This means that $M\setminus
(H_{1}^{S}\cup_{S} H_{2}^{S})$ is homeomorphic to $\partial M\times
I$.

Suppose $P$ and $S$ are well-separated, then
$H^{P}_{1}\cup_{P}H^{P}_{2}\subset \text{closed surface $\times
I$}$. Since $F^{P}_{1}\cup F^{P}_{2}$ is incompressible in $M$, each
component of $F^{P}_{i}$ is isotopic to $F\times\bigl\{0\bigr\}$ in
$F\times I$. Hence the lemma holds. Q.E.D. \vskip 2mm

The distance between two essential simple closed curves $\alpha$ and
$\beta$ on a compact surface $P$, denoted by $d(\alpha,\beta)$, is
the smallest integer $n\geq 0$ so that there is a sequence of
essential simple closed curves
$\alpha_{0}=\alpha,\ldots,\alpha_{n}=\beta$ on $P$ such that
$\alpha_{i-1}$ is disjoint from $\alpha_{i}$ for $1\leq i\leq n$.
When $P$ is a  bicompressible surface in a 3-manifold $M$, the
distance of $P$ is $d(P)=Min\bigl\{d(\alpha,\beta)\bigr\}$, where
$\alpha$ bounds a disk in $H^{P}_{1}$ and $\beta$ bounds a disk in
$H^{P}_{2}$. See [4] and [16].

Lemma 2.3 follows from the definitions and the main result in
[15]:\vskip 3mm

{\bf Lemma 2.3.} \ (1) \ If $M=V\cup_{S}W$ is a reducible Heegaard
splitting, then $d(S)=0$.

(2) \ If $M=V\cup_{S} W$ is a weakly reducible Heegaard splitting,
then $d(S)\leq 1$.

(3) \ If $M=V\cup_{S} W$ is a non-trivial and $\partial$-reducible
Heegaard splitting, then $d(S)\leq 1$.

(4) \ If $C$ is a closed surface, and $V\cup_{S} W$ is a non-trivial
Heegaard splitting of $C\times I$, then $d(S)\leq 2$ [15]. \vskip
2mm

{\bf Lemma 2.4 ([3], [13]).} \   Let $M=V\cup_{S} W$ be a Heegaard
splitting, and $P$ be an incompressible surface in $M$. Then either
$P$ can be isotoped to be disjoint from $S$ or $d(S)\leq 2-\chi(P)$.
\vskip 3mm

{\bf Lemma 2.5 ([16]).} \ Let $P$ and $Q$ be bicompressible but
weakly incompressible connected closed separating surfaces in a
3-manifold $M$. Then either

(1) \ $P$ and $Q$ are well-separated, or

(2) \ $P$ and $Q$ are isotopic, or

(3) \ $d(P)\leq 2g(Q)$. \vskip 3mm

{\bf Lemma 2.6 ([16]).} \ Let $M=V\cup_{S} W$ be a Heegaard
splitting of a 3-manifold $M$. If $d(S)>2g(S)$, then $M$ has the
unique minimal Heegaard spliting $V\cup_{S}W$ up to isotopy. \vskip
3mm

\section{The proof of Theorem 1}

Let $M=M_{1}\cup_{F} M_{2}$, and $F$ be a bounded surface. Then
$M=M^{1}\cup_{P^{1}}M^{*}\cup_{P^{2}} M^{2}$, where  $M^{1}, M^{2},
P^{1}, P^{2},
\partial_{1},
\partial_{2}$ are defined in Section 1. \vskip 3mm

{\bf Lemma 3.1.} \   $g(M)\leq g(M_{1})+g(M_{2})$.

{\bf Proof.} \ Let $M_{i}=V_{i}\cup_{S_{i}} W_{i}$ be a Heegaard
splitting of $M_{i}$ such that $F\subset
\partial_{i}\subset \partial_{-} W_{i}$. Now let $\gamma_{i}$ be a
unknotted arc in $W_{i}$ such that $\partial_{1} \gamma_{i}
\subset\partial_{+} W_{i}$,
$\partial_{2}\gamma_{1}=\partial_{2}\gamma_{2}\subset intF$. Let
$N(\gamma_{1}\cup\gamma_{2})$ be a regular neighborhood of
$\gamma_{1}\cup\gamma_{2}$ in $W_{1}\cup_{F} W_{2}$. Let
$V=V_{1}\cup N(\gamma_{1}\cup\gamma_{2})\cup V_{2}$, and $W$ be the
closure of $(W_{1}\cup_{F} W_{2})\setminus
N(\gamma_{1}\cup\gamma_{2})$. Then $V\cup_{S} W$ is a Heegaard
splitting of $M$, where $S=\partial_{+} V=\partial_{+} W$. Note that
$g(S)=g(S_{1})+g(S_{2})$. Hence the lemma holds. See also [17].
Q.E.D. \vskip 2mm

{\bf Lemma 3.2.} \  If $M_{i}$ has a Heegaard splitting
$V_{i}\cup_{S_{i}} W_{i}$ with $d(S_{i})> 2(g(S_{1})+g(S_{2}))$,
$i=1,2$, then any minimal Heegaard splitting of $M$ is irreducible
and weakly reducible.

{\bf Proof.} \ Since $d(S_{i})>2(g(S_{1})+g(S_{2}))$ for $i=1,2$,
$M_i$ is irreducible and $F$ is an essential surface in $M_i$ by
Lemma 2.3 (1) and (2). Then it follows that $M=M_1\cup _F M_2$ is
irreducible and $F$ is an essential surface in $M$. Furthermore,
$M_{i}$ is not homeomorphic to an I-bundle of a closed surface by
Lemma 2.3 (4).

Let $M=V\cup_{S} W$ be a minimal Heegaard splitting of $M$. Since
$M$ is  irreducible and $V\cup_{S} W$ is minimal, $V\cup_{S} W$ is
irreducible. By Lemma 2.6, $g(M_{i})=g(S_{i})$. By Lemma 3.1,
$g(S)\leq g(S_{1})+g(S_{2})$.

Now suppose $V\cup_{S} W$ is strongly irreducible. Since $S_{i}$ is
separating in $M$,  $S_{i}$ is bicompressible but weakly
incompressible in $M$. By Lemma 2.5, either $S$ and $S_{1}$ are
well-separated, or $S$ and $S_{1}$ are isotopic, or $d(S_{1})\leq
2g(S)$. Since $S_1$ is a Heegaard surface of $M_1$, we have that
$H^{S_{1}}_{1}\cup_{S_{1}}H^{S_{1}}_{2}$ is homeomorphic to $M_{1}$,
which not a product as we just proved. By Lemma 2.2, $S$ and $S_{1}$
are not well-separated. Since $S$ is a Heegaard surface of $M$ and
$S_{1}$ is a Heegaard surface of $M_{1}$,  $S$ is not isotopic to
$S_{1}$. Since $d(S_{1})>2(g(S_{1})+g(S_{2}))$ and $g(S)\leq
g(S_{1})+g(S_{2})$, $d(S_{1})>2g(S)$, a contradiction. Q.E.D. \vskip
3mm

{\bf The proof of Theorem 1.} \ Under the assumptions of Theorem 1,
by Lemma 3.2, any minimal Heegaard splitting of $M$ is irreducible
and weakly reducible. Let $M=V\cup_{S} W$ be a minimal Heegaard
splitting, then $V\cup _{S} W$ has a thin position as
$$V\cup_{S}W=(V_{1}^{'}\cup_{S_{1}^{'}}
W_{1}^{'})\cup_{F_{1}}\ldots\cup_{F_{n-1}}(V_{n}^{'}\cup_{S_{n}^{'}}
W_{n}^{'}) \ \ \ (*)$$ where $n\geq 2$, and each component of
$F_{1},\ldots, F_{n-1}$ is an incompressible closed surface in $M$,
and each $V_{i}^{'}\cup_{S_{i}^{'}} W_{i}^{'}$ consists of a active
component which is a non-trivial strong irreducible Heegaard
splitting and possible some product components each of which is a
trivial Heegaard splitting of an I-bundle of a closed surface. See
[14]. \vskip 3mm

Since $d(S_{i})>2(g(S_{1})+g(S_{2}))$ for $i=1,2$, $M$ is
irreducible and $F$ is an essential surface in $M$ by Lemma 2.3.
Furthermore, $M_{i}$ is not homeomorphic to an I-bundle of a closed
surface.

{\bf Claim 1.} \  $F_{i}$, $1\leq i\leq n-1$, can be isotoped so
that $F_{i}\subset M^{*}$.

{\bf Proof.} \ Suppose that $F_{i}\cap (M^{1}\cup
M^{2})\neq\emptyset$ for some $1\leq i\leq n-1$. We may assume that
$F_{i}\cap M^{1}\neq\emptyset$, and $F_{i}\cap M^{1}$ is
incompressible moreover. Note that $g(F_{i})\leq g(S)\leq
g(S_{1})+g(S_{2})$.
 Now $\chi(F_{i}\cap M^{1})\geq \chi(F_{i}\cap M_{1})\geq \chi(S)= 2-g(S) \ge
2-2(g(S_{1})+g(S_{2}))$, we have $d(S_1) > 2- \chi(F_i\cap M^{1})$.
By Lemma 2.4, $F_{i}$ can be isotoped to be disjoint from $S_{1}$.
Hence $F_{i}\cap M^{1}$ lies in one of $V_{1}$ and $W_{1}$ which
contains $F$.  It follows $F_i$ can be further isotoped  to be
disjoint from $M^{1}$. Q.E.D. (Claim 1) \vskip 3mm

 {\bf Claim 2.} \ There exists a component of $\bigcup_{1\leq i\leq n-1} F_{i}$ isotopic to $P^{1}$ (resp.
 $P^{2}$).

{\bf Proof.} \ Suppose that each component of $\bigcup_{1\leq i\leq
n-1} F_{i}$ is not isotopic to $P^{1}$. By Claim 1, $M^{1}\subset
M_{i}^{'}=V_{i}^{'}\cup_{S_{i}^{'}} W_{i}^{'}$ for some $1\leq i\leq
n$.

If $S_1$ is contained in the product component of
$M_{i}^{'}=V_{i}^{'}\cup_{S_{i}^{'}} W_{i}^{'}$, by the same
argument in the proof of Lemma 2.2,
$H^{S_{1}}_{1}\cup_{S_{1}}H^{S_{1}}_{2}$ is a Heegaard splitting of
an I-bundle of a closed surface. That is $M_1$ is an I-bundle of a
closed surface, it is a contradiction.

Now suppose $S_1$ is contained in the active component of
$M_{i}^{'}=V_{i}^{'}\cup_{S_{i}^{'}} W_{i}^{'}$. We also use
$M_{i}^{'}=V_{i}^{'}\cup_{S_{i}^{'}} W_{i}^{'}$ to denote the active
component. Now $S_{i}^{'}$ and $S_{1}$ are both weakly
incompressible in $M_{i}^{'}$ and
$H^{S_{1}}_{1}\cup_{S_{1}}H^{S_{1}}_{2}$ is not product. By Lemma
2.2, $S_{i}^{'}$ and $S_{1}$ are not well separated. By Lemma 3.1,
$g(S_{i}^{'})\leq g(S)\leq g(M_{1})+g(M_{2})$. By Lemma 2.6,
$g(M_{1})+g(M_{2})=g(S_{1})+g(S_{2})$. Hence
$2g(S_{i}^{'})<d(S_{1})$.

By Lemma 2.5, $S_{i}^{'}$ and $S_{1}$ are isotopic in $M_{i}^{'}$.

It follows $S_1$ is a strongly irreducible Heegaard surface of
$M_i'$. Then by Lemma 2.2 $M_i'\setminus
(H^{S_{1}}_{1}\cup_{S_{1}}H^{S_{1}}_{2})$ is $\partial M_i'\times
I$. By Lemma 2.2 we may assume that
$H^{S_{1}}_{1}\cup_{S_{1}}H^{S_{1}}_{2}=M^{1}$. Hence one component
of $\partial M_{i}^{'}$ is isotopic to $P^{1}$.

Similarly, there exists a component of $\bigcup_{1\leq i\leq n-1}
F_{i}$ isotopic $P^{2}$. Q.E.D. (Claim 2) \vskip 3mm

By Claim 2, Theorem 1 holds. Q.E.D. \vskip 3mm

\section{The proof of Theorem 2}


{\bf Lemma 4.1.} \ Let $N=P_{1}\times I\cup_{F} P_{2}\times I$ be
the surface sum of $P_{1}\times I$ and $P_{2}\times I$ along $F$,
where $P_{1}$ and $P_{2}$ are orientable closed surfaces, and $F$ is
a connected bounded surface in both $P_{1}\times\bigl\{0\bigr\}$ and
$P_{2}\times\bigl\{0\bigr\}$.

(1) \ If both $P_{1}\times\bigl\{0\bigr\}\setminus F$ and
$P_{2}\times\bigl\{0\bigr\}\setminus F$ are connected, then
$g(N)=Min\bigl\{g(P_{1})+g(P_{2}), \alpha\bigr\}$, where
$$\alpha=1/2(2\chi(F)+2-\chi(P_{1})-\chi(P_{2}))+Min\bigl\{g(P_{1}),
g(P_{2})\bigr\}.$$

(2) \ If $F$ is an annulus, then $g(N)=g(P_{1})+g(P_{2})$.

{\bf Proof.} \ We first prove (1).

Since both $P_{1}\times\bigl\{0\bigr\} \setminus F$ and
$P_{2}\times\bigl\{0\bigr\}\setminus F$ are connected, $N$ contains
three boundary components $P_{1}\times\bigl\{1\bigr\}$,
$P_{2}\times\bigl\{1\bigr\}$, and
$P^{*}=(P_{1}\times\bigl\{0\bigr\}\setminus F)\cup
(P_{2}\times\bigl\{0\bigr\}\setminus F)$. See Figure 4.1. Hence

$$g(N)\geq
Min\bigl\{g(P_{1})+g(P_{2}), g(P_{1})+g(P^{*}),
g(P_{2})+g(P^{*})\bigr\}.$$

\begin{center}
\psfrag{a}{$P^*$} \psfrag{c}{$P_1\times I$} \psfrag{d}{$P_2\times
I$} \psfrag{f}{$F$} \psfrag{i}{} \psfrag{j}{} \psfrag{g}{$P_1\times
1$} \psfrag{h}{$P_2\times 1$}
\includegraphics[height=2.5in]{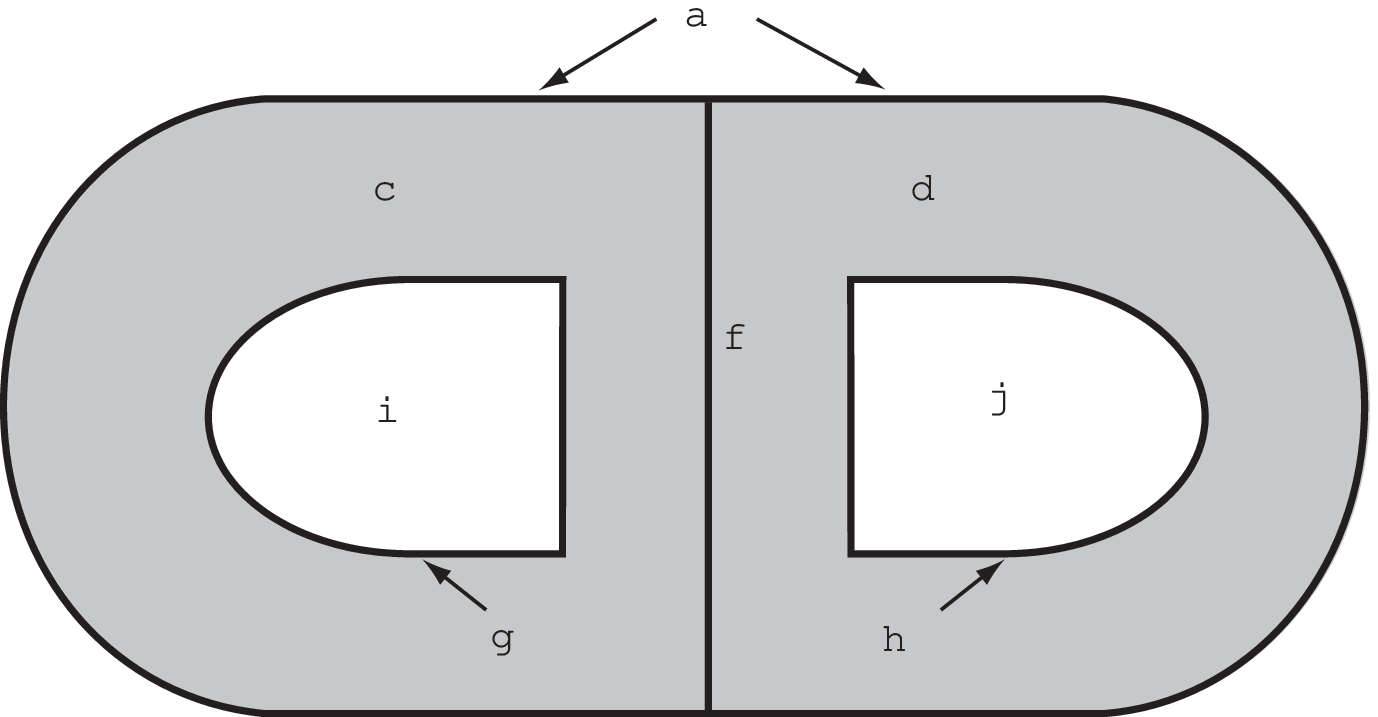}
\centerline{Figure 4.1}
\end{center}

It is easy to see that $N$ is homeomorphic to both $P_{1}\times
I\cup_{P_{1}\times\bigl\{0\bigr\}\setminus F} P^{*}\times I$ and
$P^{*}\times I\cup_{P_{2}\times\bigl\{0\bigr\} \setminus F}
P_{2}\times I$. See Figure 4.2.

\begin{center}
\psfrag{a}{$P^*\times I$}\psfrag{b}{$P_1\times I$} \psfrag{f}{$F$}
\includegraphics[height=2in]{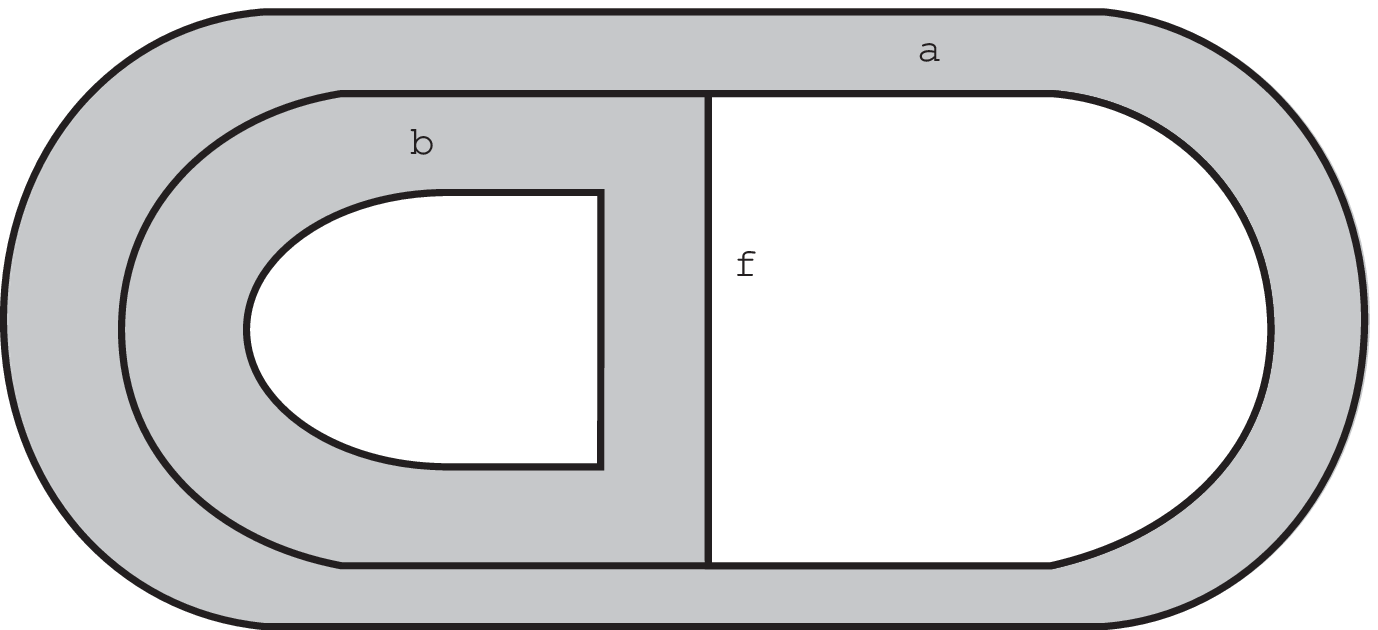}
\centerline{Figure 4.2}
\end{center}

By Lemma 3.1, $$g(N)\leq Min\bigl\{g(P_{1})+g(P_{2}),
g(P_{1})+g(P^{*}), g(P_{2})+g(P^{*})\bigr\}.$$ Now
$$g(N)=Min\bigl\{g(P_{1})+g(P_{2}),
g(P_{1})+g(P^{*}), g(P_{2})+g(P^{*})\bigr\}.$$

Since $P^{*}=(P_{1}\times\bigl\{0\bigr\}\setminus F)\cup
(P_{2}\times\bigl\{0\bigr\}\setminus F)$,
$$2-2g(P^{*})=\chi(P^{*})=\chi(P_{1})+\chi(P_{2})-2\chi(F),$$ and
$$g(P^{*})=1/2(2\chi(F)+2-\chi(P_{1})-\chi(P_{2})).$$
Hence (1) holds.

Now we prove (2).

Suppose now that $F$ is an annulus.  Now there are three cases:

Case 1. \  Both $P_{1}\times\bigl\{0\bigr\}\setminus F$ and
$P_{2}\times\bigl\{0\bigr\}\setminus F$ are connected.

Now $N$ contains three boundary components
$P_{1}\times\bigl\{1\bigr\}$, $P_{2}\times\bigl\{1\bigr\}$ and
$P^{*}=(P_{1}\times\bigl\{0\bigr\}\setminus F)\cup
(P_{2}\times\bigl\{0\bigr\}\setminus F)$. Since $F$ is an annulus,
$g(P^{*})\geq g(P_{1}), g(P_{2})$. By the argument in (1), (2)
holds.

Case 2. One of  $P_{1}\times\bigl\{0\bigr\}\setminus F$ and
$P_{2}\times\bigl\{0\bigr\}\setminus F$ is connected while the other
is non-connected.

The argument is the same with the one in Case 1.

Case 3. \  Both $P_{1}\times\bigl\{0\bigr\}\setminus F$ and
$P_{2}\times\bigl\{0\bigr\}\setminus F$ are non-connected.

Now we denote by $F_{i}^{1}$ and $F_{i}^{2}$ the two components of
$P_{i}\times\bigl\{0\bigr\}\setminus F$. We may assume that
$\partial F_{1}^{j}=\partial F_{2}^{j}$. Then $N$ contains four
boundary components $P_{1}\times\bigl\{1\bigr\}$,
$P_{2}\times\bigl\{1\bigr\}$, $F^{1}=F_{1}^{1}\cup F_{2}^{1}$ and
$F^{2}=F_{1}^{2}\cup F_{2}^{2}$. In this case,
$g(F^{1})+g(F^{2})=g(P_{1})+g(P_{2})$. Hence (2) holds. Q.E.D.
\vskip 3mm

{\bf The proof of Theorem 2.} \ Recalling the definitions of $M^{i},
M^{*}$ defined in Section 1. Since $\partial_{i}$ is separating in
$M$ for $i=1, 2$, and $M^{i}$ is homeomorphic to $M_{i}$ for
$i=1,2$, by Theorem 1,
$$g(M)=g(M_{1})+g(M_{2})+g(M^{*})-g(\partial_{1})-g(\partial_{2}).$$
By Lemma 4.1, $g(M^{*})=Min\bigl\{g(\partial_{1})+g(\partial_{2}),
\alpha\bigr\}$, where
$$\alpha=1/2(2\chi(F)+2-\chi(\partial_{1})-\chi(\partial_{2}))+Min\bigl\{g(\partial_{1}),
g(\partial_{2})\bigr\}.$$ Hence $g(M)=Min\bigl\{g(M_{1})+g(M_{2}),
\alpha\bigr\}$, where  $$\alpha=
g(M_{1})+g(M_{2})+1/2(2\chi(F)+2-\chi(\partial _{1})-\chi(\partial
_{2}))-Max\bigl\{g(\partial _{1}), g(\partial _{2})\bigl\}.$$

It is easy to see that  $g(M)=g(M_{1})+g(M_{2})$ if and only if
$\chi(F)\geq 1/2Max\bigl\{\chi(\partial_{1}),
\chi(\partial_{2})\bigr\}.$ Q.E.D.\vskip 3mm

{\bf Corollary 3.} \ The proof follows immediately from Theorem 1
and Lemma 4.1 (2). Q.E.D.\vskip 3mm

\vskip 10mm

{\bf Reference}\vskip 5mm

[1] \ D. Bachman, S. Schleimer and E. Sedgwick, Sweepouts of
amalgamated 3-manifolds, Algebr. Geom. Topol. 6(2006) 171-194.

[2] \ A. Casson and C. McA Gordon, Reducing Heegaard splittings,
Topology Appl. 27(1987) 275-283.

[3] \ K. Hartshorn, Heegaard splittings of Haken manifolds have
bounded distance, Pacific J. Math. 204(2002) 61-75.

[4] \ J. Hempel, 3-manifolds as viewed from the curve complex,
Topology 40(2001) 631-657.

[5] \ T. Kobayashi, A construction of arbitrarily high degeneration
of tunnel numbers of knots under connected sum, J. Knot Theory
Ramifications 3(1994) 179-186.

[6] \ T. Kobayashi and R. Qiu, The amalgamation of high distance
Heegaard splittings is always efficient, Math. Ann. 341(2008)
707-715.

[7] \ T. Kobayashi, R. Qiu, Y. Rieck and S. Wang, Separating
incompressible surfaces and stabilizations of Heegaard splittings,
Math. Proc. Cambridge Philos. Soc. 137(2004) 633-643.

[8] \ T. Kobayashi and Y. Rieck, Heegaard genus of the connected sum
of $m$-small knots, Commu. Anal. Geom. 14(2006), 1037-1077.

[9] \ M. Lackenby, The Heegaard genus of amalgamated 3-manifolds,
Geom. Dedicata 109(2004) 139-145.

[10] \ Tao Li, On the Heegaard splittings of amalgamated
3-manifolds, Geom. Topol. Monographs 12(2007) 157-190.

[11] \ K. Morimoto, On the super additivity of tunnel number of
knots, Math. Ann. 317(2000), 489-508.

[12] \ R. Qiu, K. Du, J. Ma and M. Zhang, Distance and the Heegaard
genera of annular 3-manifolds, Preprint.

[13] \ M. Scharlemann, Proximity in the curve complex: boundary
reduction and bicompressible surfaces, Pacific J. Math. 228(2006),
325-348.

[14] M. Scharlemann and A. Thompson, Thin position for 3-manifolds.
Geometric Topology (Haifa, 1992) 231-238, Contemp. Math., 164, Amer.
Math. Soc., Providence, RI, 1994.

[15] \ M. Scharlemann and A. Thompson, Heegaard splittings of ${\rm
Surfaces}\times I$ are standard, Math. Ann., 295(1993), 549-564.

[16] \ M. Scharlemann and M. Tomova, Alternate Heegaard genus bounds
distance, Geom. Topol. 10(2006) 593-617.

[17] \ J. Schultens, Additivity of tunnel number for small knots,
Comment. Math. Helv. 75(2000) 353-367.

[18] \ J. Schultens and R. Weidman, Destabilizing amalgamated
Heegaard splittings, Geom. Topol. Monographs 12(2007) 319-334.

[19] \ J. Souto, The Heegaard genus and distance in curve complex,
Preprint. \vskip 5mm

Department of Applied Mathematics, Dalian University of Technology,
Dalian, China

qiurf@dlut.edu.cn \vskip 5mm

School of Mathematics, Peking University, Beijing, China

wangsc@math.pku.edu.cn \vskip 5mm

Department of Applied Mathematics, Dalian University of Technology,
Dalian, China

zhangmx@dlut.edu.cn

\end{document}